# Sustainable Multi-Modal Transportation and Routing focusing on Costs and Carbon Emissions Reduction

Saba Javanpour; A. Radman; Sarow Saeedi[*]; Sina Feizi Karimabadi; Daniel A. Larson; Eric C. Jones

*Department of Industrial, Manufacturing, and Systems Engineering, University of Texas at Arlington, Texas, USA*

## Abstract

Transportation plays a critical role in supply chain networks, directly impacting cost efficiency, delivery reliability, and environmental sustainability. This study provides an enhanced optimization model for transportation planning, emphasizing environmental sustainability and cost-efficiency. An Integer Linear Programming (ILP) model was developed to minimize the total transportation costs by considering organizational and third-party vehicles' operational and rental costs while incorporating constraints on carbon emissions. The model incorporates multi-modal transportation routing and emission caps to select the optimized number of organizational and rental vehicles of different modes in each route to ensure adherence to sustainability goals. Key innovations include adding carbon emission constraints and optimizing route selection to reduce overall emissions. The model was implemented using the Gurobi solver, and numerical analysis reveals a trade-off between cost minimization and carbon footprint reduction. The results indicate that adopting tight environmental policies increases the costs by around 8% on average while more than 95% of the vehicles utilized will be rented. These insights provide actionable guidance for industries aiming to enhance both economic performance and environmental responsibility.

## 1. Introduction

Transportation involves the movement of goods, raw materials, and products between various locations within a supply chain network. It serves as a link that connects elements of the system. Effective transportation ensures the smooth flow of materials, from suppliers to manufacturers and finally to end consumers. Any disruption in transportation can cause delays, reduce overall performance, and negatively impact productivity across the entire supply chain [1].

Transportation is crucial to a company's operations and a major cost factor. Vehicle-related expenses, including fuel, maintenance, repairs, insurance, and depreciation, form a significant portion of transportation costs, especially in fleet-dependent industries[1]. These costs are influenced by various factors, including the distance traveled, vehicle type, and the age of the fleet. Efficient management of these vehicle-related costs is critical for businesses looking to optimize their overall logistics expenses.

This study focuses on a transportation system, in which the total cost of operating and renting vehicles is calculated for organization-owned and rented vehicles, respectively. The model consists of three transportation modes: truck, trailer, and van, and it considers five major Texan cities, their distances, and cost factors. To create many instances to run the model, instances and relevant factors are calculated using real-world data, with randomness added wherever applicable. To consider the sustainability development goals, CO2 emissions are also calculated using real-world data and used in the model.

---

[*] Correspondign author. Email address: sarow.saeedi@uta.edu
[1] https://fractory.com/transportation-cost-explained/ (Supply Chain Issues: Transportation Cost - Accessed in November 2024)

## 2. Literature Review

Optimal routing in supply chain networks (SCNs) involves determining the most efficient vehicle paths to minimize costs and ensure timely delivery of goods. This strategy reduces expenses related to fuel, labor, and vehicle wear and tear while meeting demand efficiently [2]. Industries like e-commerce, retail, and manufacturing rely heavily on routing optimization to enhance operational efficiency and customer satisfaction, as transportation costs can constitute up to 50% of total logistics costs [3]. Effective routing reduces costs and improves delivery reliability, enhancing customer satisfaction [4]. Additionally, it supports sustainability by cutting fuel consumption and emissions, aligning with environmental regulations and corporate responsibility goals [5]. Thus, optimal routing delivers financial and environmental benefits, making it vital in modern supply chain strategies.

The paper examines transportation modes with varying capacities, optimizing their use to prevent supply chain shortages. It enhances the base model by factoring in the environmental impact of transportation routes. Studies have linked multimodal transport and optimal routing to eco-friendly supply chains. A 2017 study found road transport accounted for 21% of the EU's $CO_2$ emission[2]. In 2021, the US transportation sector contributed 38% of the country's $CO_2$ emissions, with air and road freight leading the carbon footprint per ton-mile in 2019 [6]. Utilizing intermodal transportation helps reduce emissions [7]. Masoud and Mason [8] showed that using diverse transportation modes balances cost and delivery time through an integer linear programming model.

Further advancements involve the use of mathematical models to optimize transportation networks. Harrir and Sari-Triqui [9] developed a model that combines transport networks to minimize transport costs and $CO_2$ emissions, while Marufuzzaman and Eksioglu [10] focused on a dynamic multimodal transportation network design model to cope with supply fluctuations and natural disasters in the biofuel industry. Liu et al. [11] estimated long-haul freight emissions, showing a 28% reduction by 2050 through mode shifts. Environmental constraints are central to studies like Lam and GuLam and Gu [12], who proposed a bi-objective intermodal model with carbon restrictions, and, Ma et al. [13] , who developed a multi-objective model for cold supply chains under governmental carbon limits.

Transportation studies use methodologies like mathematical modeling, heuristic solutions, and multi-criteria decision analysis (MCDA). Models such as integer linear and mixed-integer nonlinear programming optimize multi-mode routing and logistics. Heuristic methods like simulated annealing and Benders decomposition address complex NP-hard problems [8], [10], [14]. Additionally, MCDA methods are utilized to select the most appropriate transportation networks by evaluating multiple criteria to make informed decisions [15]. Researchers have explored multi-mode transportation's role in supply chain sustainability using methods like surveys and numerical analysis [16], life cycle assessment [17], and discrete event simulation ([18]) to evaluate environmental impacts.

The paper introduces the capacitated transportation modes as a literature gap by its publication time. However, as we showed earlier, at the time of the paper's publication, there existed studies that considered transportation modes in supply chains.

## 3. Methodology
### 3.1. Methodological Overview

In this Project, we utilized the following framework to approach the problem. As shown in Figure *1*, each route i-j has a demand for different types of transportation modes, which can be satisfied by the transportation organization in different periods with rental, organization-owned vehicles, or both. Multiple modes of transportation that are rented or owned by the transportation organization have a stopping duration in each source city and some other costs in the route i-j that can be minimized. Considering the limited budget of the organization in source cities to rent multiple modes of transportation, our model determines the optimal number of rental and organizational vehicles of different modes.

---

[2] https://climate.ec.europa.eu/eu-action/transport-emissions/road-transport-reducing-co2-emissions-vehicles_en (Accessed: October 2024)

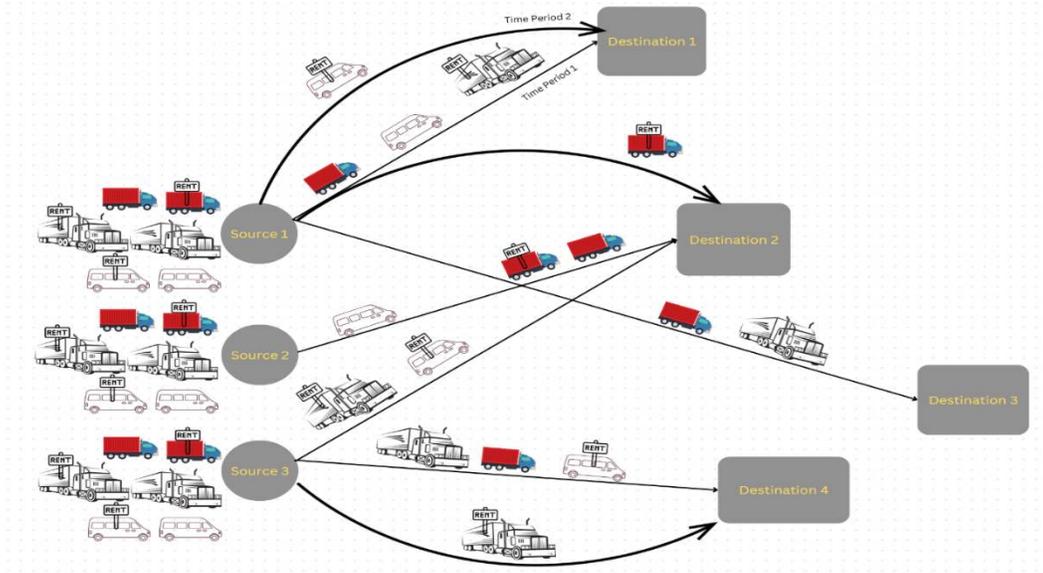

Figure 1. Multimodal Transportation network with rental and organization-owned vehicles

### 3.2. Basic Model's Components

- **Sets and Indices:**

| | |
|---|---|
| $m$ | Index of vehicle types including truck, trailer, van, etc., $m = \{1, \dots, M\}$ |
| $n$ | Index of transportation centers, $n = \{1, \dots, N\}$ |
| $i$ | Index of origins, $i = \{1, \dots, I\}$ |
| $j$ | Index of destinations, $j = \{1, \dots, J\}$ |
| $t$ | Index of the time period, $t = \{1, \dots, T\}$ |

- **Parameters:**

| | |
|---|---|
| $V_{imt}$ | maximum number of the organizational vehicle $m$ in origin $i$ at period $t$ |
| $Vr_{imt}$ | maximum number of the rental vehicle $m$ in origin $i$ at period $t$ |
| $D_{ijmt}$ | required number of vehicles of type $m$ in $i-j$ route at period $t$ |
| $CS_{imt}$ | cost associated with the stopping duration of the organizational vehicle of type $m$ in the origin $i$ at period $t$ |
| $CSr_{imt}$ | cost associated with the stopping duration of the rental vehicle of type $m$ in the origin $i$ at period $t$ |
| $CT_{ijmt}$ | cost associated with traveling from origin $i$ to destination $j$ by the organizational vehicle of type $m$ at period $t$ |
| $CTr_{ijmt}$ | cost associated with traveling from origin $i$ to destination $j$ by the rental vehicle of type $m$ at period $t$ |
| $TC$ | Required budget for renting and operational costs |
| $OPR_{imt}$ | operational cost of a vehicle of type $m$ in origin $i$ at period $t$ |
| $Rent_{imt}$ | cost associated with renting a vehicle of type m in origin $i$ at period $t$ |

- **Decision Variables:**

| | |
|---|---|
| $x_{ijmt}$ | Number of organizational vehicle type $m$ travels from $i$ to $j$ at time $t$ |
| $xr_{ijmt}$ | Number of rental vehicle type $m$ travels from $i$ to $j$ at time $t$ |
| $y_{imt}$ | Number of organizational vehicle type $m$ is empty in origin $i$ at period $t$ |
| $yr_{imt}$ | Number of rental vehicle type $m$ is empty in origin $i$ at period $t$ |
| $q_{imt}$ | Number of organizational vehicle type $m$ is in service in origin $i$ at period $t$ |
| $qr_{imt}$ | Number of rental vehicle type $m$ is in service in origin $i$ at period $t$ |

- **Objective Function**

The objective function minimizes the total cost of transportation by organizational and rental vehicles and minimizes the stopping duration of vehicles in unloading locations.

$$Min\ z = \sum_{i=1}^{I}\sum_{j=1}^{J}\sum_{m=1}^{M}\sum_{t=1}^{T} CT_{ijmt} \times x_{ijmt} + \sum_{i=1}^{I}\sum_{j=1}^{J}\sum_{m=1}^{M}\sum_{t=1}^{T} CTr_{ijmt} \times xr_{ijmt}$$
$$+ \sum_{i=1}^{I}\sum_{m=1}^{M}\sum_{t=1}^{T} CS_{imt} \times y_{imt} + \sum_{i=1}^{I}\sum_{m=1}^{M}\sum_{t=1}^{T} CSr_{imt} \times yr_{imt} \quad (1)$$
$$+ \sum_{i=1}^{I}\sum_{m=1}^{M}\sum_{t=1}^{T} OPR_{imt} \times q_{imt} + \sum_{i=1}^{I}\sum_{m=1}^{M}\sum_{t=1}^{T} Rent_{imt} \times qr_{imt}$$

- **Constraints:**

$$\sum_{j=1}^{J} x_{ijmt} + y_{imt} = V_{imt} \qquad \forall\ i,m,t \quad (2)$$

$$\sum_{j=1}^{J} xr_{ijmt} + yr_{imt} = Vr_{imt} \qquad \forall\ i,m,t \quad (3)$$

$$\sum_{i=1}^{I}\sum_{m=1}^{M} x_{ijmt} + \sum_{i=1}^{I}\sum_{m=1}^{M} xr_{ijmt} \geq \sum_{i=1}^{I}\sum_{m=1}^{M} D_{ijmt} \qquad \forall\ j,t \quad (4)$$

$$\sum_{i=1}^{I}\sum_{m=1}^{M}\sum_{t=1}^{T} OPR_{imt} \times q_{imt} + \sum_{i=1}^{I}\sum_{m=1}^{M}\sum_{t=1}^{T} Rent_{imt} \times qr_{imt} \leq TC \quad (5)$$

$$\sum_{j=1}^{J} x_{ijmt} \leq q_{imt} \qquad \forall\ i,m,t \quad (6)$$

$$\sum_{j=1}^{J} xr_{ijmt} \leq qr_{imt} \qquad \forall\ i,m,t \quad (7)$$

### 3.3. Modified Model

- **New Parameters**

$E_m$      Emission factor (kg $CO_2$ per km) for organizational vehicle type $m$
$Er_m$      Emission factor (kg $CO_2$ per km) for rental vehicle type $m$
$D_{ij}$      Distance between origin $i$ and destination $j$
$Emission\ Cap$      Maximum allowable emissions for the entire transportation network

- **Modified Sustainability Constraints:**

The following constraints limit the overall carbon emissions to a certain threshold ensuring the transportation emissions remain within acceptable environmental limits:

$$\sum_{i,j,m,t} (E_m \times D_{ij} \times x_{ijmt} + Er_m \times D_{ij} \times xr_{jmt}) \leq Emission\ Cap \quad (8)$$

In Eq. (8), emissions for organizational vehicles equal to $E_m \times D_{ij} \times x_{ijmt}$, and emissions for rental vehicles equal to $Er_m \times D_{ij} \times xr_{ijmt}$. By adding Eq. (8), the model ensures that the transportation network follows the emission cap policy.

### 3.4. Data Collection

We utilized the following data collection methods in this paper:

*Primary Data Extraction from Literature and Standards*: Parameters like freight costs are derived from logistics standards and accessible databases. These sources offer valuable baseline figures for costs, applicable to both organizational and rental vehicles[3] [19].

*Model Data Generation*: For variables such as maximum vehicle capacity, organizational budget, and vehicle operational cost, data is generated through calculated estimates and historical values [3].

## 4. Results and Discussion

To compare the base model with our enhanced model, we utilized the numerical data which is uploaded in the following path: https://github.com/SarowSaeedi/Transportation-Routing.git.

We solved the models in the Gurobi solver, and the results for both the base and the enhanced models for a problem with the size of 2*3*2*2 ($I*J*M*T$) are demonstrated in Table *1***Error! Reference source not found.**. The base model's objective function value reflects a purely cost-centric optimization. In contrast, the enhanced model, which integrates sustainability considerations by limiting overall carbon emissions, has a slightly higher objective function value.

Table 1. Comparison of the base model and the enhanced model for a problem size of 2*3*2*2

| | Base Model | | | | Enhanced Model | | |
|---|---|---|---|---|---|---|---|
| Objective Function = 16,792.69 | x[1,1,1,1]=1 | x[2,3,2,2]=1 | yr[2,2,1]=3 | Objective Function = 17,169.77 | x[1,1,1,1]=1 | xr[2,1,1,2]=1 | yr[2,2,2]=3 |
| | x[1,1,1,2]=2 | xr[2,1,1,2]=1 | yr[2,2,2]=2 | | x[1,1,1,2]=2 | xr[2,1,2,1]=1 | q[1,1,1]=4 |
| | x[1,1,2,1]=1 | xr[2,1,2,2]=1 | q[1,1,1]=4 | | x[1,2,1,1]=2 | xr[2,3,1,1]=2 | q[1,1,2]=4 |
| | x[1,2,1,1]=2 | xr[2,3,1,1]=2 | q[1,1,2]=4 | | x[1,2,1,2]=2 | xr[2,3,1,2]=1 | q[1,2,1]=1 |
| | x[1,2,1,2]=2 | xr[2,3,1,2]=1 | q[1,2,1]=2 | | x[1,2,2,1]=1 | y[1,2,1]=3 | q[1,2,2]=2 |
| | x[1,2,2,1]=1 | y[1,2,1]=2 | q[1,2,2]=2 | | x[1,2,2,2]=2 | y[2,1,1]=2 | q[2,2,1]=3 |
| | x[1,2,2,2]=2 | y[2,1,1]=2 | q[2,2,1]=3 | | x[1,3,1,1]=1 | y[2,1,2]=2 | q[2,2,2]=3 |
| | x[1,3,1,1]=1 | y[2,1,2]=2 | q[2,2,2]=3 | | x[2,1,2,1]=1 | yr[1,1,1]=4 | qr[1,2,2]=1 |
| | x[2,1,2,1]=1 | yr[1,1,1]=4 | qr[2,1,1]=2 | | x[2,1,2,2]=2 | yr[1,1,2]=4 | qr[2,1,1]=2 |
| | x[2,1,2,2]=1 | yr[1,1,2]=4 | qr[2,1,2]=2 | | x[2,3,2,1]=2 | yr[1,2,1]=4 | qr[2,1,2]=2 |
| | x[2,2,2,2]=1 | yr[1,2,1]=4 | y[2,1,1]=2 | | x[2,3,2,2]=1 | yr[1,2,2]=1 | qr[2,2,1]=1 |
| | x[2,3,2,1]=2 | yr[1,2,2]=2 | qr[2,2,2]=1 | | xr[1,2,2,2]=1 | yr[2,2,1]=2 | y[2,1,2]=2 |

### 4.1. Sensitivity Analysis on *Emission Cap*:

To observe the impact of the *Emission Cap* values, we ran the enhanced model 15 times with different values of the *Emission Cap*. The results of the sensitivity analysis, shown in Figure *2* (a), demonstrate that as the emission cap becomes more stringent (lower values), the objective function value for the developed model increases notably. By running the enhanced model 15 times, we can see an average of 8% increase in the costs compared to the base model.

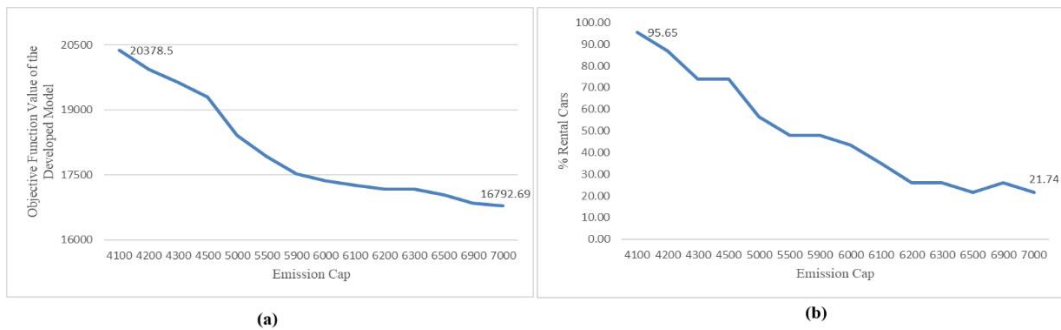

Figure 2. Effect of Emission Cap Values on the Objective Function (a), and Rate of Rental Cars (b).

---

[3] https://www.bts.gov/ (Freight Facts and Figures – Accessed: November 2024)
  https://www.fmcsa.dot.gov/ (Hours of Service (HOS) | FMCSA - Department of Transportation – Accessed November 2024).

This indicates that achieving stricter environmental goals involves higher transportation costs, reflecting the economic trade-offs required for carbon emission reductions. On the other hand, Figure *2* (b) demonstrates that by reducing the Emission Cap, the model tends to select more rental vehicles, as rental vehicles are newer than their organizational counterparts and consequently produce fewer emissions.

### 4.2. Conclusion

This study demonstrates the effect of applying eco-friendly policies in multimodal a transportation network. We enhanced a mathematical model incorporating practical considerations such as emission caps and rental car usage. The enhanced model's costs are around 8% higher than the base model which does not commit to the environmental policies, providing a trade-off for decision-making in the transportation network. The sensitivity analysis further highlights the effectiveness of the enhanced model in balancing economic and environmental trade-offs. As emission caps increase, the objective function value systematically decreases, showing a clear trade-off between profitability and sustainability. Simultaneously, the percentage of rental cars decreased significantly, indicating the enhanced model's efficiency in reducing dependency on high-emission vehicles under stricter environmental regulations. These results validate our enhanced model's ability to integrate sustainability constraints into production and distribution planning, offering a robust framework for achieving operational efficiency while adhering to environmental goals. Future research may expand on this approach by incorporating the sink nodes' demands as parameters and transportation types as decision variables.